\documentclass{ajour}
\usepackage{epsfig,amsmath,amssymb}

\newenvironment{itemise}
 {\begin{itemize}}{\end{itemize}}

\newcommand{\al}{{\alpha}}
\newcommand{\be}{\beta}
\newcommand{\ga}{\gamma}
\newcommand{\del}{\delta}

\newcommand{\ep}{\varepsilon}
\newcommand{\la}{\lambda}

\newcommand{\Q}{\mathbb{Q}}

\newcommand{\ds}{\oplus}
\newcommand{\bigds}{\bigoplus}
\newcommand{\tens}{\otimes}

\newcommand{\iso}{\cong}

\newcommand{\sdp}{\rtimes}

\newcommand{\bmat}{\begin{pmatrix}}
\newcommand{\emat}{\end{pmatrix}}

\newcommand{\gl}[2]{\text{\rm GL}_{#1}\!\left(#2\right)}

\newcommand{\glnk}{\gl{n}{k}}
\newcommand{\agl}[2]{\text{\rm AGL}_{#1}\!\left(#2\right)}
\newcommand{\aglnk}{\agl{n}{k}}
\newcommand{\malg}[2]{\text{\rm M}_{#1}\!\left(#2\right)}
\newcommand{\mpr}[3]{{{#1}[{#2}]_{#3}}}

\DeclareMathOperator{\rad}{rad}

\DeclareMathOperator{\End}{End}

\DeclareMathOperator{\diag}{diag}

\begin{document}

%


\authorrunninghead{Scott H. Murray}
\titlerunninghead{Conjugacy classes in parabolic subgroups}



\title{Conjugacy classes in maximal parabolic subgroups\\ of general linear groups}
\author{Scott H. Murray} 
\date{\today}
\affil{University of Chicago}
\email{murray@math.uchicago.edu}

\abstract
{We compute conjugacy classes in maximal parabolic
subgroups of the general linear group. This computation proceeds by
reducing to a ``matrix problem''.  Such problems involve finding
normal forms for matrices under a specified set of row and column
operations.  We solve the relevant matrix problem in small dimensional
cases.  This gives us all conjugacy classes in maximal
parabolic subgroups over a perfect field when one of the two blocks
has dimension less than 6.  In particular, this includes every maximal 
parabolic subgroup of $\glnk$ for $n<12$ and $k$ a perfect field.  If
our field is finite of size $q$, we also show that the number of 
conjugacy classes, and so the number of characters, of these groups is 
a polynomial in $q$ with integral coefficients.}


\begin{article}

\section{Introduction}
A great deal of progress has been made recently towards describing the
representation theory of reductive
algebraic groups.  For example, the study of representations of finite
reductive groups was greatly advanced by the work of Deligne and Lusztig
\cite{mr52:14076} and has been an active field of research.
Conjugacy classes in reductive groups have been investigated by
Springer and Steinberg \cite{MR42:3091}.  In comparison, little is known for
solvable algebraic groups \cite{MR99b:20011}. Even less is known
about  groups which are neither reductive nor solvable.   

The parabolic subgroups of the
general linear group are among the simplest such ``mixed'' groups.
Each is a semidirect product of the unipotent radical (which is a solvable 
normal subgroup) with a Levi complement (which is a
reductive group).  Representations of the Levi complement can be
inflated to the maximal parabolic---this is 
vital to the inductive step of classifications of representations of the general
linear group.  Drozd \cite{mr94d:22016} generalized these inflated
representations to a much larger class of mixed groups and showed they
are Zariski dense in the set of  irreducible representations.
Almost nothing is known about the other representations of the
maximal parabolics.  

For parabolic subgroups of a reductive group, the
conjugacy classes contained in the unipotent radical were first
investigated by Richardson, R{\"o}hrle, and Steinberg
\cite{mr93j:20092}.  A series of papers on this subject have been
written by Hille, J\"urgens, Popov, and R\"ohrle 
\cite{MR1669178,MR1692308,MR1665003,MR99f:14063,MR1617826,MR97c:20070,MR99b:20081,MR99e:20060,MR1669615}.
Some of their results use matrix problems similar to the ones
discussed here.

In this paper, we describe conjugacy classes in maximal parabolic
subgroups.  This can be 
considered as a  step towards a better understanding of the
representation theory of parabolic subgroups of reductive algebraic
groups.  A maximal parabolic subgroup 
$$G = P^{(m,n)}=
\bmat\gl{m}{k}&\malg{m,n}{k}\\0&\gl{n}{k}\emat$$
has {\em unipotent
radical\/}\index{unipotent radical}  
$$U=\bmat I_m&\malg{m,n}{k}\\0&I_n\emat$$
and {\em Levi complement\/}\index{Levi complement}
$$L=\bmat\gl{m}{k}&0\\0&\gl{n}{k}\emat.$$
Multiplication in the unipotent radical is given by
$$\bmat I_m & v \\ 0 & I_n \emat\bmat I_m & w \\ 0 & I_n \emat=\bmat I_m & v+w \\ 0 & I_n \emat.$$
Hence $U$ can be identified with $M_{m,n}(k)$, the additive group of
$m\times n$ matrices
The Levi subgroup acts on the unipotent radical in the natural manner:
$$\bmat A&0\\0&B\emat \bmat 1&v\\0&1\emat \bmat A^{-1}&0\\0&B^{-1}\emat
=\bmat  1&AvB^{-1}\\0&1\emat.$$

The following lemma describes the conjugacy classes in
a semidirect product with abelian normal subgroup.  This is analogous to the
description of the characters proved by Clifford theory 
\cite[section~11B]{mr90k:20001}. 
\begin{lemma}\label{L-ccsdp}
Let the group $G$ be a semidirect product $U\sdp L$ with $U$ abelian.
Then, for every $h$ in $L$, the conjugacy classes in $G$ intersecting
$Uh$ are in one-to-one correspondence with the orbits of $C_L(h)$ on
$C^U(h)=U/[U,h]$. 
\end{lemma}
\noindent This is proved by taking $u,v\in U$ and $h,k\in L$ and then rearranging
$(vk)(uh)(vk)^{-1}$ to be a product of an element of $U$ with an element
of $L$.

This lemma provides us with a procedure for finding the conjugacy
classes.  In Section~\ref{S-genjcf} we describe the generalized Jordan
normal form, which gives us a set of conjugacy class
representatives for $L$.  Then, for every
such representative $h$, we compute the  centralizer 
$C_L(h)$ in Sections~\ref{S-cent} and \ref{S-gens}, and the
cocentralizer $C^U(h)$ in Section~\ref{S-cocent}.   Note that
Sections~\ref{S-genjcf}, \ref{S-cent} and \ref{S-cocent} each have two
subsections: in the first we consider matrices with 
rational eigenvalues, in the second we show that for 
an irrational separable eigenvalue we get essentially the same thing,
but over 
the extension of $k$ with the eigenvalue adjoined.  If you are only
interested in algebraic closed fields, you need
only read the first subsection in each section.

Finding orbits of the centralizer on the cocentralizer
turns out to be a ``matrix problem''. 
Such problems involve finding normal forms for matrices under a
specified set of row and column operations.  They have been
extensively studied by the Kiev school 
founded by Nazarova and Ro\u\i ter \cite{mr49:4877}.  A good
reference on matrix problems and their applications to representations 
of algebras is \cite{mr98e:16014}.
There is a classification of matrix
problems:  finite type problems have finitely many
orbits whose representatives can be 
independent of the field;  while for infinite type, the number of
orbits depends on the field and is infinite whenever the
field is.   Infinite type problems can further be divided into tame
type, where an explicit solution is known; and wild type, which
reduces to the classical unsolved problem of finding a normal form for a
pair of noncommuting matrices.  The Brauer-Thrall conjecture, proved
in \cite{MR54:360}, shows
that every matrix problem is either finite, tame or wild.

In Section~\ref{S-solmp},  we  show that the  matrix
problems associated with $P^{(m,n)}$ are all finite type if, and only
if,  $m$ or $n$ is less than 6.   In particular, this includes every
maximal  parabolic subgroup of $\glnk$ for $n<12$ and $k$ a perfect
field.    If our field is finite of size $q$, it follows from our
proof that the number of  conjugacy classes, and so the number of
characters, of these groups is  a polynomial in $q$ with integral
coefficients. Finally in Section~\ref{S-ccagl} we recompute the
conjugacy classes of the affine general linear groups.

\section{Notation}\label{S-nota}
This section explains some of the notational conventions used in
this paper.  
We use standard group theoretic notation as in \cite{mr96m:20001}.
We use the equality symbol to denote  natural isomorphism as well as
strict equality.

Throughout this paper $k$\index{$k$} is a field.  Many of our
results work for perfect fields or separable field extensions only.
We denote the  
set of $m\times n$  
matrices over a $k$-algebra $R$ by $\malg{m,n}{R}$\index{$\malg{m,n}{R}$} and write $\malg{n}{R}$\index{$\malg{n}{R}$}
for $\malg{n,n}{R}$.  We frequently consider the {\em multiplicative group}\index{multiplicative group} $R^\times$
of a $k$-algebra\index{$R^\times$}; for example, the {\em general linear group}\index{general linear group}
$\glnk$\index{$\glnk$} is just $\malg{n}{k}^\times$.   The $n$-dimensional space of
column vectors is $k^n$\index{$k^n$}.  We write $A^t$\index{$A^t$} for the {\em transpose}\index{transpose}
of the matrix $A$\index{$A$}.  We  also apply this operation to sets of
matrices, for example, $(k^n)^t$\index{$(k^n)^t$} is the space of row vectors over $k$.
The algebra direct sum of $n$ copies of $k$ is written $k^{\ds n}$ and
is identified with the algebra of {\em diagonal matrices}\index{diagonal matrices} in $\malg{n}{k}$.

A {\em composition}\index{composition} of $n$ is a sequence
$\la=(\la_1,\la_2,\dots,\la_s)$\index{$\la=(\la_1,\la_2,\dots,\la_s)$} of natural numbers whose sum is
$n=|\la|$. We call $\la_i$ the $i$th part of $\la$.
A {\em partition}\index{partition} is a composition whose parts are in decreasing order.
Often we find it convenient to
write the partition $\la$ in the form
$(r^{l_r},\dots,2^{l_2},1^{l_1})$\index{$\la=(r^{l_r},\dots,2^{l_2},1^{l_1})$} where $l_j$ is 
the number of times the part $j$ occurs.

We  consider algebras or groups consisting of matrices with
different kinds of entries in different positions.  These are denoted by
matrices whose entries are the appropriate sets of possible entries.
A matrix 
whose entries are also matrices is called a {\em block matrix}, and is
identified in the obvious manner with a matrix of larger dimension.
For example, we have defined the elements of $P^{(n,m)}$ as $2\times2$
matrices of matrices, but we generally consider them as 
$(n+m)\times(n+m)$ matrices over $k$.

We define the Jordan block
$$J_n(A)=\bmat
A &0 & 0&\cdots&0\\
I_d&A&0&\cdots&0\\
0&I_d&A&\ddots&\vdots\\
\vdots&\vdots&\ddots&\ddots&0\\
0&0&\cdots&I_d&A\emat$$
where $A$ is a $d\times d$ matrix and  $n$ is the number of times $A$ appears.
We also write $J_\la(A) = \bigds_{i=1}^s J_{\la_i}(A)$ for a
composition $\la$

\section{Generalized Jordan normal
  form\index{generalized Jordan normal form}}\label{S-genjcf}
We need a set of conjugacy class representatives
for the Levi complement $L$, so that we can apply Lemma~\ref{L-ccsdp}.
Since this group is just 
$\gl{m}{k}\ds\glnk$, it suffices to give representatives of the
similarity classes of invertible matrices.  For our
purposes, these representatives should be rational\index{rational} (ie.\ defined over
$k$) but also as close to diagonal as possible.  The generalized
Jordan normal form has these properties.  A proof that any
matrix over a perfect field is similar to a matrix in this form
can be found in 
\cite{mr29:3477}.   

Corresponding to every direct sum decomposition of a matrix, there is a
decomposition of the underlying vector space into a direct sum of
subspaces invariant under that matrix.
An element of $\malg{n}{k}$ which is not similar to a direct sum of
smaller square matrices over $k$ is called 
{\em indecomposable}\index{indecomposable matrix}.  Every square matrix is
a direct sum of indecomposables. 

Fix an $n\times n$ matrix $A$.  Given $p$, a monic irreducible
polynomial over $k$, the subspace
 $V_p=\left\{ v\in k^n : p(A)^mv=0 \text{ for some natural number $m$}\right\}$
is easily seen to be $A$-invariant.  If this subspace is nonzero,
we say $p$ is a {\em generalized eigenvalue\/}\index{generalized
  eigenvalue} of $A$ and $V_p$ is the corresponding {\em generalized
eigenspace\/}\index{generalized eigenspace}.  In fact, $k^n$ is a
direct sum of the generalized eigenspaces, so we get a corresponding decomposition
of $A$.  In particular, an indecomposable matrix has a unique
generalized eigenvalue whose generalized eigenspace is the entire
underlying vector space.

\subsection*{Rational case}
Suppose every  generalized eigenvalue of $A$ is of the form
$p(t)=t-\al$, for some $\al$ in $k$.  Then each such $\al$ is also an
eigenvalue. So $A$ is
similar to $\bigds_\al A_\al$ where $A_\al$ has unique eigenvalue $\al$.
The indecomposables with eigenvalue $\al$ are just the Jordan blocks
$J_m(\al)$.  Hence $A_\al$ is similar to  $J_{\la_\al}(\al)$
for some partition $\la_\al$ and $A$ is similar to 
$\bigds_\al J_{\la_\al}(\al)$
where $\al$ runs over the eigenvalues of $A$ and
$n=\sum_{\al}\left|\la_\al\right|$. Of course, this is just the ordinary Jordan
normal form.

\subsection*{Irrational case}
Now suppose the generalized eigenvalues of $A$  are arbitrary monic,
separable, irreducible polynomials. 

Let $p$ be such a polynomial and write $K=k(\al)$, where $\al$ is a
root of $p$ in the algebraic closure 
$\bar{k}$.  Then $K$ is a separable field extension of $k$ and multiplication by
$\al$ induces a $k$-linear 
transformation $K\to K$.  The matrix of this transformation with respect
to the $k$-basis $\{1,\al,\dots,\al^{d-1}\}$ is just the companion matrix
of $p$, denoted $C_p$.
This is an indecomposable matrix with generalized eigenvalue $p$.
More generally, every indecomposable matrix with generalized eigenvalue $p$
is similar to the matrix $J_m(C_p)$ for some natural number $m$.

So the matrix $A$ is similar to $\bigds_{p} J_{\la_p}(C_p)$,
where $p$ runs over the generalized eigenvalues of $A$ and the $\la_p$ are
partitions with $n=\sum_{p}\left|\la_p\right|\deg(p)$.  This is the
{\em generalized Jordan normal form}\index{generalized Jordan normal form} of $A$.

\section{Centralizers in general linear groups}\label{S-cent}
We describe the centralizers in general
linear groups.  Our description is explicit provided that
all the generalized eigenvalues of our matrix are separable over
$k$.  We first compute the centralizer in $M=\malg{n}{k}$, then use
this to find the centralizer in $G=\glnk$. 

Let  $A=\bigds_{p} J_{\la_p}(C_p)$ be an element of $G$ in generalized
Jordan normal form.  Corresponding to this decomposition of $A$ is a
decomposition of  $k^n$ into a direct sum of  generalized 
eigenspaces $V_p$.  Further,
a matrix $B$ that centralizes $A$ also centralizes $p(A)$, and so $V_p$ is
invariant under $B$.  Hence $C_M(A) = \bigds_p
C_{\malg{n_p}{k}}(J_{\la_p}(C_p))$ where $n_p=|\la_p|\deg(p)$.
We may now assume, without loss of generality, that $A$ has a single
generalized eigenvalue $p$ with corresponding partition $\la=\la_p$.

\subsection*{Rational case}
First consider  $p(t)=t-\al$, ie.\ the 
matrix has rational eigenvalue $\al$.  Now  $A=J_\la(\al)=\al I_n +
J_\la(0)$ and $\al I_n$ is in the center of $M$, so
$C_M(A) = C_M(J_\la(0))$.    
We write $J_\la$ for $J_\la(0)$ and $k[x]_n$ for $k[x]/(x^n)$.  
With $\la=(n)$ there is an isomorphism
$C_M(J_n)\to k[x]_n$ taking $J_n$ to $x$. 
The natural  action of
$C_M(J_n)$  on $k^n$ is equivalent to the regular action of
$k[x]_n$ on itself.  We  generalize this to an arbitrary partition.

\begin{figure}
\centerline{\raisebox{1in}{$\left(\rule{0in}{1.1in}\right.$}
\psfig{figure=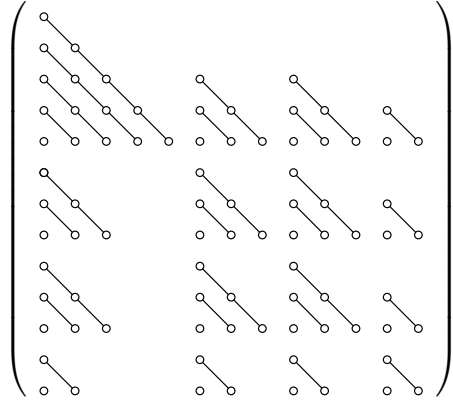}
\raisebox{1in}{$\left.\rule{0in}{1.1in}\right)$}}
\caption{An element of the centralizer}\label{F-cent}
\end{figure}

Take a matrix $B$ centralizing $J_\la$ and write it in block form
  $$B= \left( \begin{matrix}
                B_{11}&\cdots& B_{1s}\\
                \vdots&\ddots&\vdots\\
                B_{s1}&\cdots&B_{ss}
              \end{matrix} \right),$$
where $B_{ij}$ is a $\la_i \times \la_j$ matrix. Then $BJ_\la = J_\la B$
implies $B_{ij}J_{\la_j} = J_{\la_i}B_{ij}$
for all $i$ and $j$.
If we write $B_{ij} = (b_{l,m})$, this becomes $b_{l,m+1} = b_{l-1,m}$ and
$b_{1\la_j}=b_{l,\la_j}=b_{\la_i,m}=0$ for $l=2,\dots,\la_i$ and
$m=1,\dots,\la_j-1$.  Hence $B_{ij}$ is
 $$
 \bmat
  \begin{matrix}
   b_0&0&\cdots&0\\
   b_1&b_0&\ddots&\vdots\\
   \vdots &\ddots &\ddots&0\\
   b_{\la_i}&\cdots&b_1&b_0\end{matrix}& \text{\Large0}
 \emat \text{ or }
 \bmat
   \text{\Large0} \\ 
   \begin{matrix}
   b_{\la_j-\la_i}&0&\cdots&0\\
   b_{\la_j-\la_i+1}&b_{\la_j-\la_i}&\ddots&\vdots\\
   \vdots &\ddots &\ddots&0\\
   b_{\la_j}&\cdots&b_{\la_j-\la_i+1}&b_{\la_j-\la_i}\end{matrix} \emat $$
for $\la_i\le\la_j$ or $\la_i\ge\la_j$ respectively.
The appearance of the full matrix $B$ is illustrated by Figure~\ref{F-cent}
for $\la=(5,3,3,2)$, where each line represents entries which are
equal and blank spaces represent zero entries.
Define $X_{c\times d}^a$ to be the $c\times d$ matrix whose $(i,j)$-entry
is 1 if $i=j+a$ and 0 otherwise.  Then $B_{ij}$ can be
written as  $\sum_{a=0}^{\la_i} b_aX_{\la_i\times\la_j}^a$ or
$\sum_{a=\la_j-\la_i}^{\la_j} b_aX_{\la_i\times\la_j}^a$ respectively.  
It is easily checked that 
$X_{c\times d}^aX_{d\times e}^b=X_{c\times e}^{a+b}$ for any
nonnegative integers $a,b,c,d,e$.  This identity gives us an algebra
homomorphism 
$$\bmat
  k[x]& x^{\la_1-\la_2}k[x]& \cdots & x^{\la_1-\la_s}k[x]\\
  k[x]& k[x]               & \cdots & x^{\la_2-\la_s}k[x]\\
\vdots& \vdots             & \ddots& \vdots\\
  k[x]& k[x]               & \cdots & k[x]
\emat\to C_M(A)$$
which takes $x^a$ in the $(i,j)$-entry to $X_{\la_i\times\la_j}^a$ in the $(i,j)$-block.  
The matrices $X_{c\times d}^a$ are linearly independent for
$a=0,\dots,\min(c,d)$ and zero for $a>\min(c,d)$.  So this
homomorphism is surjective and its kernel is 
$$\bmat
  (x^{\la_1})& (x^{\la_1})& \cdots  & (x^{\la_1})\\
  (x^{\la_2})& (x^{\la_2})& \cdots  & (x^{\la_2})\\
  \vdots     & \vdots     & \ddots  & \vdots\\
  (x^{\la_s})& (x^{\la_s})& \cdots  & (x^{\la_s})
\emat.$$
Hence $C_M(J_\la)$ is isomorphic to the quotient algebra\index{$\mpr{k}{x}{\la}$}
$$\mpr{k}{x}{\la}=\bmat
  k[x]_{\la_1}& x^{\la_1-\la_2}k[x]_{\la_1}& \cdots & x^{\la_1-\la_s}k[x]_{\la_1}\\
  k[x]_{\la_2}& k[x]_{\la_2}               & \cdots & x^{\la_2-\la_s}k[x]_{\la_2}\\
  \vdots      & \vdots                     & \ddots & \vdots\\
  k[x]_{\la_s}& k[x]_{\la_s}               & \cdots & k[x]_{\la_s}
\emat.$$

Next we consider the natural action of this algebra.  The left action of
$C_M(J_\la)$ on $k^n$ is identical to its action
on its own first column.  This is easily seen
to be isomorphic to the action of $\mpr{k}{x}{\la}$ on its first
column which is
$$\mpr{k}{x}{\la\times-}=\bmat k[x]_{\la_1}\\\vdots\\k[x]_{\la_s}\emat.$$
We also need the right action on row vectors.  There is an isomorphism
$D:\mpr{k}{x}{\la}\to\mpr{k}{x}{\la}^t$\index{$D:\mpr{k}{x}{\la}\to\mpr{k}{x}{\la}^t$} given by
$$\bmat
  a_{11}& x^{\la_1-\la_2}a_{12}& \cdots & x^{\la_1-\la_s}a_{1s}\\
  a_{21}& a_{22}                 & \cdots & x^{\la_2-\la_s}a_{2s}\\
  \vdots      & \vdots                     & \ddots & \vdots\\
  a_{s1}& a_{s2}               & \cdots & a_{ss}
\emat\mapsto
\bmat
  a_{11}& a_{12}& \cdots & a_{1s}\\
  x^{\la_1-\la_2}a_{21}& a_{22}                 & \cdots & a_{2s}\\
  \vdots      & \vdots                     & \ddots & \vdots\\
  x^{\la_1-\la_s}a_{s1}& x^{\la_2-\la_s}a_{s2}               & \cdots & a_{ss}
\emat.$$
We consider $\mpr{k}{x}{\la}$ to act naturally on
$\mpr{k}{x}{-\times\la}=(\mpr{k}{x}{\la\times-})^t$\index{$\mpr{k}{x}{-\times\la}$}
via this isomorphism.

We now turn to the multiplicative group of $\mpr{k}{x}{\la}$, which
is isomorphic to $C_G(A)=C_M(A)^\times$.  
Now  $B\in\mpr{k}{x}{\la}$ can be written
  $$B = B_0 + (xI)B_1 + (xI)^2B_2 +\cdots + (xI)^{\la_s-1}B_{\la_s-1},$$
where the entries of each $B_i$ are in $k$ and $xI$ is the matrix with
$x$ in each diagonal entry and zero elsewhere. 
But $(xI)^{\la_s}=0$, so $xI$
is nilpotent and $B$ is invertible exactly when $B_0$ is.  Writing
$\la$ as $(r^{l_r},\dots,2^{l_2},1^{l_1})$, we see that $B_0$
is block lower triangular of the form
  $$B_0= \bmat
                B_{rr}& \cdots& 0\\
                \vdots&\ddots&\vdots\\
                B_{1r}&\cdots &B_{11}
              \emat,$$   
where $B_{ij}$ is an $l_i \times l_j$ matrix.  Hence $B$ is
invertible if, and only if, all the matrices $B_{ii}$ are invertible, and we have
described $\mpr{k}{x}{\la}^\times$\index{$\mpr{k}{x}{\la}^\times$}.

\subsection*{Irrational case}
Finally we tackle the case where $p$ is  nonlinear and separable. 
Let $\al=\al_1,\al_2,\dots,\al_d$ be the distinct roots of $p$ in 
the algebraic closure $\bar{k}$.  We take $A$ to be $J_\la(C_p)$.
Recall that $C_p$ is the $k$-matrix of multiplication by $\al$ on
$K=k(\al)$ with respect to the basis $\{1,\al,\dots,\al^{d-1}\}$.
Hence $C_{\malg{d}{k}}(C_p)$ is isomorphic to the centralizer of $\al$
in $\End_k(K)$, which is just $\End_K(K)=K$.  On the other hand, the
Jordan normal form of $C_p$ is $D_p=\diag(\al_1,\dots,\al_d)$, so
$tD_pt^{-1}=C_p$ for some $t$ in $GL_d(\bar{k})$.  The centralizer of
$D_p$ in $\malg{d}{\bar{k}}$ is the algebra of diagonal matrices
$\bar{k}^{\ds d}$.  Hence
$$C_{\malg{d}{\bar{k}}}(C_p) = t(C_{\malg{d}{\bar{k}}}(D_p))t^{-1} = t(\bar{k}^{\ds d})t^{-1},$$
so $C_{\malg{d}{k}}(C_p)\iso K$ is the set of elements of
$t(\bar{k}^{\ds d})t^{-1}$ defined over $k$.

Turning now to $J_\la(C_p)$, we see it is conjugated to $J_\la(D_p)$ by
$t^{\ds n}$.  Further,  $J_\la(D_p)$ is conjugated to
$J_\la(\al_1)\ds\cdots\ds J_\la(\al_d)$ by the obvious permutation of
basis elements.  From the rational case, we know that the centralizer
in $\malg{n}{\bar{k}}$ of this last matrix is isomorphic to
$\mpr{\bar{k}}{x}{\la}^{\ds d}$.  Reversing these
conjugations we get $C_{\malg{n}{\bar{k}}}(J_\la(D_p))\iso
\mpr{(\bar{k}^{\ds d})}{x}{\la}$ and  $C_{\malg{n}{\bar{k}}}(J_\la(C_p))\iso
\mpr{\left(t(\bar{k}^{\ds d})t^{-1}\right)}{x}{\la}$.  So $C_M(J_\la(C))$ is just the
set of elements of this algebra defined over $k$, which is
 $\mpr{K}{x}{\la}$.  The natural action and
multiplicative group can now be computed as in the rational case.

An example should make this process clearer.  Suppose $k=\Q$,
$p(t)=t^2-3$, and $\la=(2)$.  Then $K=\Q(\surd3)$,
$\al_1=\surd3$ and $\al_2=-\surd3$. Hence $J_\la(C_p)$, $J_\la(D_p)$, and
$J_\la(\al_1)\ds J_\la(\al_2)$ are
$$
\left(\begin{array}{cc|cc} 0 & 1 & 1 & 0 \\
               -3& 0 & 0 & 1 \\\hline
               0 & 0 & 0 & 1 \\
               0 & 0 & -3& 0 \end{array}\right),
\left(\begin{array}{cc|cc}\surd3&0& 1 & 0 \\
              0&-\surd3& 0 & 1 \\\hline
               0 & 0 & \surd3&0 \\
               0 & 0 & 0&-\surd3 \end{array}\right) \text{, and }
\left(\begin{array}{cc|cc}\surd3&1& 0 & 0 \\
              0&\surd3& 0 & 0 \\\hline
               0 & 0 &-\surd3&1 \\
               0 & 0 & 0&-\surd3 \end{array}\right)$$
 respectively.  The centralizers of $J_\la(\al_1)\ds J_\la(\al_2)$
and $J_\la(D_p)$ consist of matrices of the form
$$\left(\begin{array}{cc|cc} a & b & 0 & 0 \\
               0 & a & 0 & 0 \\\hline
               0 & 0 & c & d \\
               0 & 0 & 0 & c \end{array}\right) \text{and}
\left(\begin{array}{cc|cc} a & 0 & b & 0 \\
               0 & c & 0 & d \\\hline
               0 & 0 & a & 0 \\
               0 & 0 & 0 & c \end{array}\right)$$ respectively.
Finally $C_M(J_\la(C_p))\iso K[x]_\la=K[x]_2$.

\section{Generators for the centralizers}\label{S-gens}
In order to find the orbits of the centralizers on the
cocentralizers, we need  a generating set for the
centralizers.  The generators we use are analogous to the elementary
matrices of linear algebra---thus finding 
the orbits  becomes a matrix problem.  

Using the notation of the previous
section, $C_G(A)$  is isomorphic to the group
$\mpr{K}{x}{\la}^\times$, which we write in block form as
$$\bmat
  \malg{l_r}{\mpr{K}{x}{r}}^\times &\cdots&\malg{l_rl_2}{x^{r-2}\mpr{K}{x}{r}} &\malg{l_rl_1}{x^{r-1}\mpr{K}{x}{r}}\\
  \vdots                            &\ddots&\vdots                               &\vdots\\
  \malg{l_2l_r}{\mpr{K}{x}{2}}     &\cdots&\malg{l_2}{\mpr{K}{x}{2}}^\times      &\malg{l_2l_1}{x\mpr{K}{x}{2}}\\
  \malg{l_1l_r}{K}                  &\cdots&\malg{l_1l_2}{K}                     &\malg{l_1}{K}^\times
\emat$$
with $\malg{l}{K[x]_r}^\times = \gl{l}{K}+\malg{l}{xK[x]_r}$.

Define the
following matrices in $\mpr{K}{x}{\la}^\times$: 
\begin{itemise}
\item $M_{i,a}(l)$, for $i=1,\dots,r$, $a\in
  {k[x]_{l_i}}^\times$, and $l=1,\dots,l_i$: diagonal entries all 1, except for the
  $(l,l)$-entry in the $(i,i)$-block which is equal to $a$;
  off-diagonal entries all 0.
\item $E_i(l,m)$, for $i=1,\dots,r$, and $l,m =1,\dots,l_i$: 
  entries all 1 and  off diagonal entries all 0, except in the
  $(i,i)$-block where the $(l,l)$ and $(m,m)$-entries are 0 and the $(l,m)$ 
  and $(m,l)$-entries are 1.
\item $A_{i\le j,a}(l,m)$, for $i, j=1,\dots,r$ with $i\le j$, $a\in k[x]_{i}$,
  $l=1,\dots,l_i$ and $m=1,\dots,l_j$: diagonal entries all 1;  off diagonal entries
  all zero, except in the $(i,j)$-block where the $(l,m)$-entry is 
  $x^{i-j}a$.
\item $A_{i\ge j,a}(l,m)$, for $i, j=1,\dots,r$ with $i\ge j$,  $a\in
  k[x]_{j}$, $l=1,\dots,l_i$ and $m  =1,\dots,l_j$: diagonal entries
  all 1;  off diagonal entries 
  all zero, except in the $(i,j)$-block where the $(l,m)$-entry is
  $a$.
\end{itemise}
Denote the $l$th row in the $i$th block by $R_{i,l}$, and the $l$th
column in the $i$th block by $C_{i,l}$.  Then these matrices act on
$\mpr{K}{x}{\la}^\times$ as in Table~\ref{T-regops}.
\begin{table}
\caption{Regular action of $\mpr{K}{x}{\la}^\times$}\label{T-regops}
\begin{tabular*}{\textwidth}{@{\extracolsep{\fill}}lll}\hline
Matrix              & Row operation & Column operation \\\hline
$M_{i,a}(l)$        & $R_{i,l}\to a\cdot R_{i,l}$ & $C_{i,l}\to C_{i,l}\cdot a$\\
$E_i(l,m)$          & $R_{i,l}\leftrightarrow R_{i,m}$ & $C_{i,l}\leftrightarrow C_{i,m}$\\
$A_{i\le j,a}(l,m)$ & $R_{i,l}\to R_{i,l}+ax^{i-j}\cdot R_{j,m}$ & $C_{j,m}\to C_{j,m}+C_{i,l}\cdot ax^{i-j}$\\
$A_{i\ge j,a}(l,m)$ & $R_{i,l}\to R_{i,l}+a\cdot R_{j,m}$ & $C_{j,m}\to C_{j,m}+C_{i,l}\cdot a$\\\hline
\end{tabular*}
\end{table}

In order to prove that these matrices generate $\mpr{K}{x}{\la}^\times$, it
suffices to show that we can reduce any matrix in this group to
the identity using these row and column operations. We proceed by
induction on the number of parts of $\la$ (ie.\ the dimension of our
matrices).  The result is clear if $\la$ has one part.  Now take a matrix 
 $$B = B_0 + (xI)B_1 + (xI)^2B_2 +\cdots + (xI)^{r-1}B_{r-1}$$
in $\mpr{K}{x}{\la}^\times$ and write 
  $$B_0= \bmat
                B_{rr}& \cdots& 0\\
                \vdots&\ddots&\vdots\\
                B_{1r}&\cdots &B_{11}
              \emat$$   
as is Section~\ref{S-cent}.
Since $B_0$ is invertible, we know $B_{rr}$ is invertible, and so we can
use row and column operations within the first block to get the
(1,1)-entry of $B$ to be one.  Now, by adding multiples of the top
row to the other rows, we can make every other entry in the first
column zero.  We can also add multiples of the left-most column of $B$ to
a column in the $i$th block, as long as we also multiply by $x^{r-i}$.
This is not a problem  since an entry in the $i$th block of the top row must be
a multiple of $x^{r-i}$ anyway.  We can now ignore the first row and column and
reduce the rest of the matrix to the identity by induction.  Hence we
are done. 

In subsequent sections we study the natural action
of these matrices on the spaces $\mpr{K}{x}{\la\times-}$ and
$\mpr{K}{x}{-\times\la}$.  These actions are  slightly different from the
regular action, because the right action is via the isomorphism $D$ of
Section~\ref{S-cent}.  The row and column operations for the natural
action are in Table~\ref{T-natops}.
\begin{table}
\caption{Natural action of $\mpr{K}{x}{\la}^\times$}\label{T-natops}
\begin{tabular*}{\textwidth}{@{\extracolsep{\fill}}lll}\hline
Matrix              & Row operation & Column operation \\\hline
$M_{i,a}(l)$        & $R_{i,l}\to a\cdot R_{i,l}$ & $C_{i,l}\to C_{i,l}\cdot a$\\
$E_i(l,m)$          & $R_{i,l}\leftrightarrow R_{i,m}$ & $C_{i,l}\leftrightarrow C_{i,m}$\\
$A_{i\le j,a}(l,m)$ & $R_{i,l}\to R_{i,l}+ax^{i-j}\cdot R_{j,m}$ & $C_{i,l}\to C_{i,l}+C_{j,m}\cdot ax^{i-j}$\\
$A_{i\ge j,a}(l,m)$ & $R_{i,l}\to R_{i,l}+a\cdot R_{j,m}$ & $C_{i,l}\to C_{i,l}+C_{j,m}\cdot a$\\\hline
\end{tabular*}
\end{table}

\section{Cocentralizers}\label{S-cocent}
Now that we have the centralizer and its generators in terms of
algebras, we find a similar description for the cocentralizer and
the action of the generators on it.

Let $G=P^{(m,n)}=U\sdp
L$.  The Levi complement is $L=\gl{m}{k}\ds \glnk$ and we can identify
the unipotent radical $U$ with the
additive group of $\malg{m,n}{k}$.  Note that
$\malg{m,n}{k}=k^m\tens(k^n)^t$\index{$\tens$ algebra tensor product over $k$} where we are, as always, tensoring over
$k$. Let $h=A\ds B$, where $A\in\gl{m}{k}$  
and $B\in\glnk$ 
are both in generalized Jordan normal form. 
The action of $L$ on $U$ is given by $A\ds B\cdot v=AvB^{-1}$.  We wish to
describe $C^U(h)=U/[U,h]$ as a $C_L(h)$-module over $k$.  
Using transfer of structure and the fact that
$I\ds B^{-1}$ is in the center of $C_L(h)$, this module is isomorphic
to $(I\ds B^{-1})\cdot C^U(h)= U/((1\ds B^{-1})\cdot[U,A\ds B])$. Finally,
\begin{eqnarray*}
 (I\ds B^{-1})\cdot[U,A\ds B] &=& \{(I\ds B^{-1})\cdot(v-A\ds B\cdot v):v\in U\} \\
                    &=& \{vB -Av:v\in U\}.
\end{eqnarray*}

We know that $C_L(h)$ is a direct sum of $C_{{\rm GL}_m(k)}(A)$ and
$C_{{\rm GL}_n(k)}(B)$, each of which is a direct sum of centralizers
corresponding to the generalized eigenvalues of $A$ and $B$.  So a
matrix $v$ in $(I\ds B^{-1})[U,A\ds B]$ can be broken up into blocks corresponding to these
centralizers and we can treat each block separately.  Hence we may assume,
without loss of generality, that $A$ and $B$ each have a single
generalized eigenvalue. 

Take $A=J_\mu(C_p)$ and $B=J_\nu(C_q)$ where $p$ and $q$ are monic,
separable, irreducible polynomials.  Let $K=k(\al)$ and
$K'=k(\be)$ with $\al,\be\in\bar{k}$ the roots 
of $p$ and $q$ respectively.   We identify 
$C_L(A\ds B)=C_{\gl{m}{k}}(A)\ds C_{\glnk}(B)$ with
$\mpr{K}{x}{\mu}^\times\ds\mpr{K'}{y}{\nu}^\times$.  This allows us to identify $U=k^m\tens(k^n)^t$  with 
$$\mpr{K}{x}{\mu\times-}\tens\mpr{K'}{y}{-\times\nu} = 
\bmat
K[x]_{\mu_1}\\\vdots\\K[x]_{\mu_s}\emat\tens
\bmat
K'[y]_{\nu_1}& \cdots& K'[y]_{\nu_t}\emat = 
\bmat
 R_{11}' & \cdots & R_{1t}' \\
                \vdots  & \ddots & \vdots \\
                R_{s1}' & \cdots & R_{st}'\emat $$
where
$$R_{ij}'=K[x]/(x^{\mu_i})\tens K'[y]/(y^{\nu_j})= K\tens K'[x,y]/(x^{\mu_i},y^{\nu_j}).$$
Now the action of $A$ on $k^m$ corresponds to the action of $(\al+x)I_m$ on
$\mpr{K}{x}{\mu\times-}$ and the action of $B$ on $(k^n)^t$ corresponds to $(\be+y)I_n$ on $\mpr{K'}{y}{-\times\nu}$.
Hence $(I\ds B^{-1})\cdot[U,A\ds B]$ is identified with the set of elements of the form 
$$v(\be+y)I_m-(\al+x)I_nv = (\be-\al+y-x)v$$ 
for $v\in \mpr{K}{x}{\mu\times-}\tens\mpr{K'}{y}{-\times\nu}$.  Hence the 
the $(i,j)$-entry of $C^U(h) \iso U/(I\ds B^{-1})\cdot[U,A\ds B]$ is identified with
$$R_{ij}= R_{ij}'/(\be-\al+y-x)=K\tens K'[x,y]/(x^{\mu_i},y^{\nu_j},\be-\al+y-x).$$

\subsection*{Rational case}
Suppose that $\al$ and $\be$ are both in $k$. Then $K=K'=k$ and
$$ R=R_{ij} = k[x,y]/(x^{\mu_i},y^{\nu_j},\be-\al+y-x),$$
If $\al\ne\be$, then $\rad R$ contains $x$, $y$ and $\be-\al = x-y$.
So the head of $R$, $R/\!\rad R$, maps onto $k/(\be-\al)=0$ and hence $R=0$.

So we can assume $\al = \be$.  Then $x=y$ in $R$, so $R=k[x]/(x^{\mu_i},x^{\nu_j})=k[x]_{l_{ij}}$ where $l_{ij}$
is the minimum of $\mu_i$ and $\nu_j$.  
Hence $C^U(h)$ becomes\index{$\mpr{k}{x}{\mu\times\nu}$}
$$\mpr{k}{x}{\mu\times\nu} =\bmat
  k[x]_{l_{11}} &  k[x]_{l_{12}} & \cdots & k[x]_{l_{1s}}\\ 
  k[x]_{l_{21}} &  k[x]_{l_{12}} & \cdots & k[x]_{l_{2s}}\\ 
  \vdots        & \vdots         & \ddots & \vdots\\ 
  k[x]_{l_{r1}} &  k[x]_{l_{r2}} & \cdots & k[x]_{l_{rs}}\emat,\quad l_{ij}=\min(\mu_i,\nu_j).$$
Since $x$ and $y$ are
identified, we have $C_L(h)\iso\mpr{k}{x}{\mu}^\times \ds \mpr{k}{y}{\nu}^\times=\mpr{k}{x}{\mu}^\times \ds \mpr{k}{x}{\nu}^\times$
acting on $C^U(h)\iso\mpr{k}{x}{\mu\times\nu}$.

\subsection*{Irrational case}
Now consider arbitrary monic, separable, irreducible polynomials $p$
and $q$.  Since $K'= k[u]/(q(u))$, we have 
\begin{eqnarray*}
  R =R_{ij}&=& K[x,y,u]/(q(u),x^{\mu_i},y^{\nu_j},u-\al+y-x)\\
    &=& K[x,y]/(q(x-y+\al),x^{\mu_i},y^{\nu_j}).
\end{eqnarray*}
Over the field $K$,  $q(u)=(u-\al)^\ep f(u)$, where $\ep$ is 1 or 0
depending on whether $p$ and $q$ are equal or unequal.  In either case
$f(\al)\ne0$.  So
we have $q(x-y+\al)= (x-y)^\ep f(x-y+\al)$ and $\rad((x-y)^\ep,
f(x-y+\al))=K[x,y]$ as it contains $x-y$ and so also contains
$f(\al)$, which is a unit.  Hence $((x-y)^\ep,f(x-y+\al))=K[x,y]$ and,
by the Chinese Remainder theorem \cite[Section {III}.2]{MR86j:00003}, 
$$R =  K[x,y]/((x-y)^\ep,x^{\mu_i},y^{\nu_j})\ds
K[x,y]/(f(x-y+\al),x^{\mu_i},y^{\nu_j}).$$
The second summand is trivial since its radical contains $x$ and $y$,
so its head maps onto $K/(f(\al))=0$.  The first summand is trivial for
$\ep=0$ and is $K[x]_{l_{ij}}$ for $\ep=1$.  Hence $C^U(h)$ is trivial for
$p\ne q$ and is isomorphic to $\mpr{K}{x}{\mu\times\nu}$ for $p=q$.
Once again $C_L(h)$ can be identified with $\mpr{K}{x}{\mu}^\times \ds
\mpr{K}{x}{\nu}^\times$. 
\subsection*{}
So we have reduced our problem to finding the orbits of
$\mpr{K}{x}{\mu}^\times\ds\mpr{K}{x}{\nu}^\times$ on
$\mpr{K}{x}{\mu\times\nu}$, for appropriate 
fields $K$.  Further, this action is given by the row and column
operations of Table~\ref{T-natops}.  So we have reduced to a matrix
problem, which we also denote $\mpr{K}{x}{\mu\times\nu}$.

\section{Solving the matrix problem for small dimensions}\label{S-solmp}
We solve the matrix problem
$\mpr{k}{x}{\mu\times\nu}$ described in the previous sections for an
arbitrary field $k$ and either $|\mu|$ or $|\nu|$ less than 6.  In
particular this gives us the conjugacy classes in maximal parabolics
of the general linear group of dimension less than 12 over a perfect field.  

\begin{figure}
\centerline{\psfig{figure=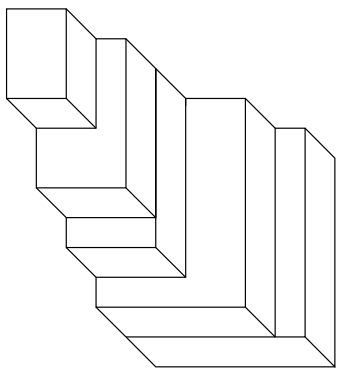}}
\caption{A three dimensional array}\label{F-3dp}
\end{figure}

Let $\mu=(r^{m_r},\dots,2^{m_2},1^{m_1})$ and
$\nu=(s^{n_s},\dots,2^{n_2},1^{n_1})$ be a pair of partitions with
$m=|\mu|$ and $n=|\nu|$.  We wish to find a normal form for matrices in
$$\mpr{k}{x}{\mu\times\nu} = \bmat
  \malg{m_rn_s}{\mpr{k}{x}{\min(r,s)}} & \cdots & \malg{m_rn_2}{\mpr{k}{x}{2}} & \malg{m_rn_1}{k}\\
  \vdots                                & \ddots & \vdots                        & \vdots\\
  \malg{m_2n_s}{\mpr{k}{x}{2}}         & \cdots & \malg{m_2n_2}{\mpr{k}{x}{2}} & \malg{m_2n_1}{k}\\
  \malg{m_1n_s}{k}         & \cdots & \malg{m_1n_2}{k} & \malg{m_1n_1}{k}\\
\emat$$
under the row and column operations of Table~\ref{T-natops}.
I find it useful to visualize such a matrix as a three dimensional
array of elements of $k$, with rows and columns as usual, and levels corresponding to the
powers of $x$.  This array is not  
rectangular since the number of levels depends on which row and column
you are in.  Figure~\ref{F-3dp} illustrates such an array for
$\mu=(6,5^2,4^2,3,2)$ and $\nu=(5^2,4,2^2,1)$.    
So, for example,
multiplying a row by $1+x^a$  takes every level
in that row and adds its entries $i$ levels higher up in the same row.
Note that to add a column to another column $i$ blocks to the left, we
also have to move $i$ levels up.  We don't have this complication when adding
to a column on the right, although we cannot add to a lower level.
Similar we need to move to a higher level when adding a row to
another row above it.

We now prove that our matrix problem can be infinite type when $m=n=6$.
\begin{theorem}\label{T-it}
The matrix problem $\mpr{k}{x}{(4,2)\times(4,2)}$ is infinite type.
\end{theorem}
\begin{proof}
Consider matrices in $\mpr{k}{x}{(4,2)\times(4,2)}$ of the form
  $$\left( \begin{matrix}
  \al x^{2}+\cdots& \be x+\cdots\\
  \ga x+ \cdots & \del +\cdots
  \end{matrix} \right)$$
for $\al,\be,\ga,\del$ in $k^\times$.
It is easily checked that every
allowable row or column operation leads to another matrix of the same
form and preserves the value of $\al\be^{-1}\ga^{-1}\del$. Hence
there are at least as many orbits as elements of $k^\times$, and the problem
is infinite type.
\end{proof}

Next we prove our main theorem, showing that all smaller matrix
problems are finite type. 
\begin{theorem}\label{T-ft}
The matrix problem $\mpr{k}{x}{\mu\times\nu}$ is finite type for
$\nu$ arbitrary and 
$\mu$ of the form $(2^{m_2},1^{m_1})$, $(r,1^{m_1})$ or
$(3,2)$.  In particular, $\mpr{k}{x}{\mu\times\nu}$ is finite type whenever $|\mu|<6$.
\end{theorem}
\begin{proof} Our basic approach is to solve the 0th level using the
permissable row and column operations, then to solve the 1st level
using only those operations which preserve the 0th level, and so on.
It is a general property of finite type matrix problems that solutions
can be found with every entry either 0 or 1.  We call positions with a
1 entry {\em pivots}\index{pivot}.  These pivots can be used to
``kill'' other positions (ie.\ make them 0 with a row or column operation).

The proof is in four cases:
\begin{enumerate}
\item 
First we consider $\mu=(2^{m_2},1^{m_1})$.  The solution 
for the 0th level is shown in the cutaway diagram of
Figure~\ref{F-mp21-arblvl0}.  
\begin{figure}\begin{center}
\label{F-mp21-arb}
\psfig{figure=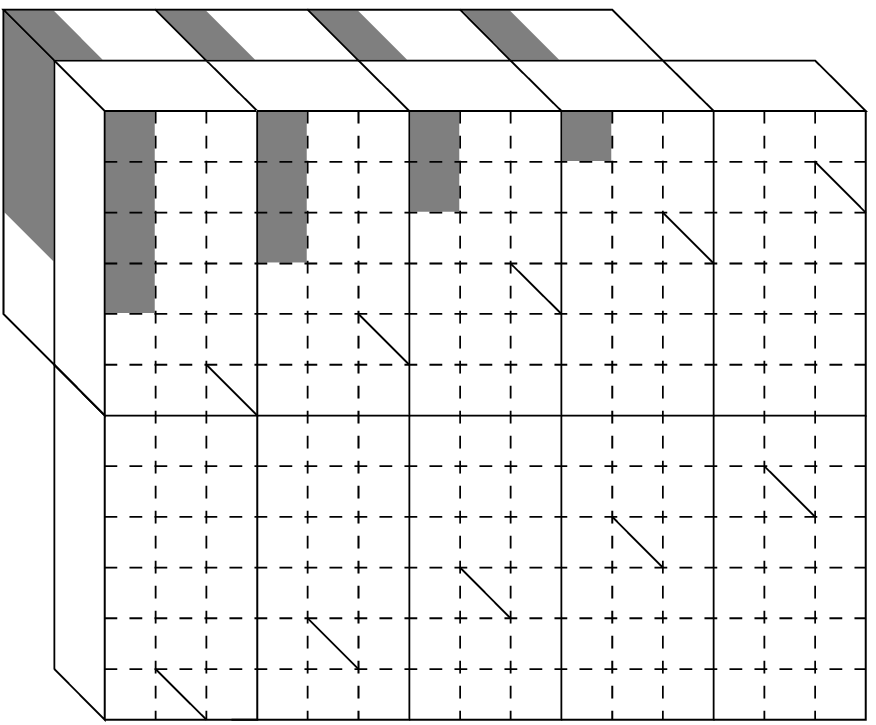}
\letteredcaption{a}{$\mu=(2^{m_2},1^{m_1})$, Level 0}\label{F-mp21-arblvl0}
\vspace{.5cm}
\psfig{figure=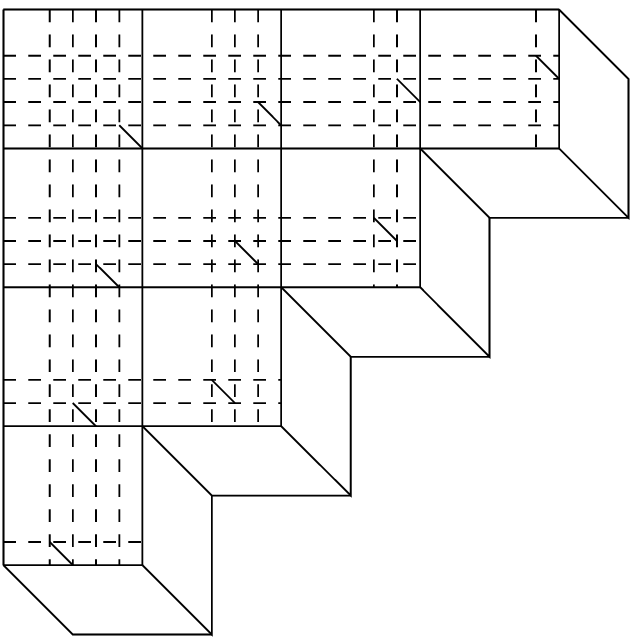}
\letteredcaption{b}{$\mu=(2^{m_2},1^{m_1})$, Level 1}\label{F-mp21-arblvl1}
\end{center}
\end{figure}
In this diagram
$\nu=(5^{m_5},\dots,1^{m_1})$ but the general
case is easily seen to be similar.  The blocks are divided by solid lines.  Each
square containing a diagonal line is an identity matrix (of course, they are not all actually
the same size).  Now we can use the pivots in the 0th
level to kill everything in the 1st level, except for the shaded blocks

The shaded blocks of the 1st level are redrawn in
Figure~\ref{F-mp21-arblvl1}.  We can add
columns to blocks on the left, but not on the right.  Also we can add
rows to blocks below but not above, because, when adding to a block
below, the damage done by one pivot can be repaired by a column operation
from another pivot.  This level can now be solved as shown.

\item 
Next we consider $\mu=(r)$, which is shown in Figure~\ref{F-mpr-arb}.
\begin{figure}
\centerline{\psfig{figure=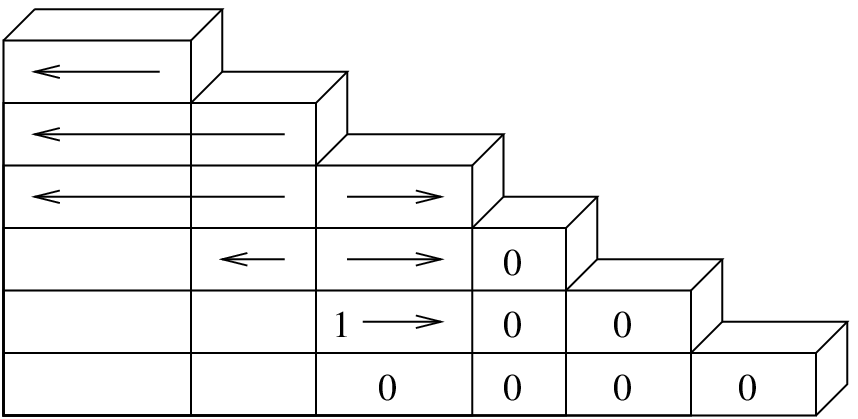}}
\caption{$\mu=(r)$}\label{F-mpr-arb}
\end{figure}
Find the first nonzero block starting in the bottom left as shown. We
can use row operations to put a pivot at the left hand end of this block and then
kill the other entries indicated by the arrows.  Now ignore all the
entries marked with an arrow or a zero, and repeat the same
process with what remains.

\item 
The case $\mu=(r,1)$ is illustrated in Figure~\ref{F-mpr1-arb}.
\begin{figure}
\centerline{\psfig{figure=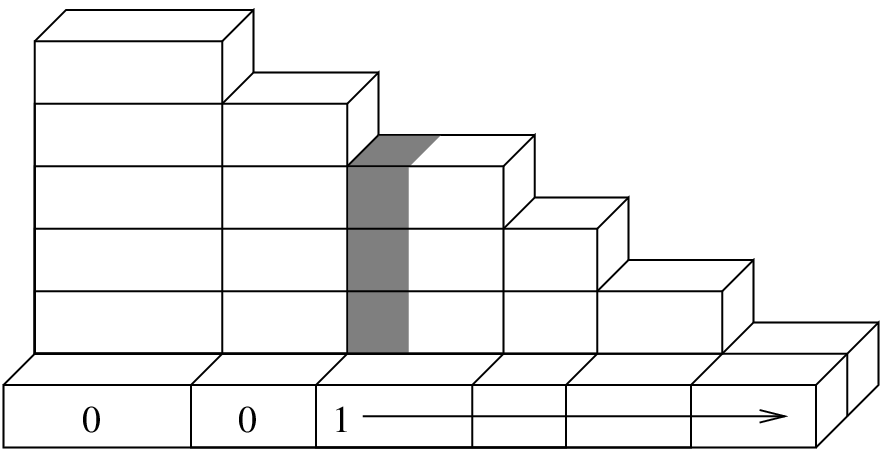}}
\caption{$\mu=(r,1)$}\label{F-mpr1-arb}
\end{figure}
We start by solving the 2nd row: find the first nonzero block, make a
pivot in that block and kill the rest of the row.
Then,  ignoring the shaded part, we solve the rest of
the 1st row as with
$\mu=(r)$.  We can now use a column multiplication to ensure that the
shaded positions contains a single 1, followed by a row multiplication to
repair any damage this does to the pivot in the second row.

The solution for $\mu=(r,1^{m_1})$ is easily seen to be similar,
except that there can be more than one shaded column.

\item
Finally we turn to {$\mu=(3,2)$}.
\begin{figure}\begin{center}
\psfig{figure=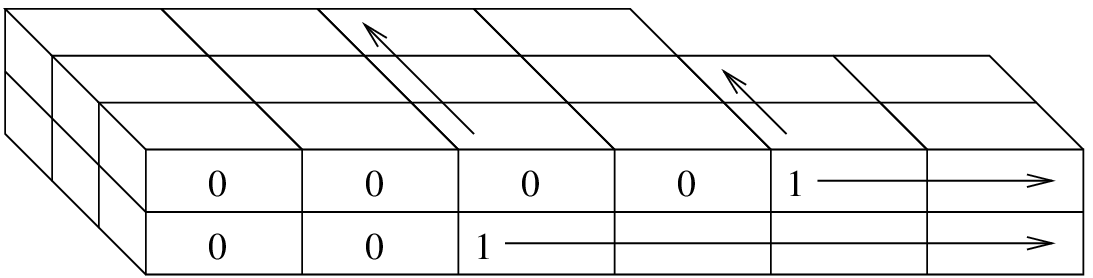}
\letteredcaption{a}{$\mu=(3,2)$, Level 0}\label{F-mp32-arblvl0}
\vspace{.5cm}
\psfig{figure=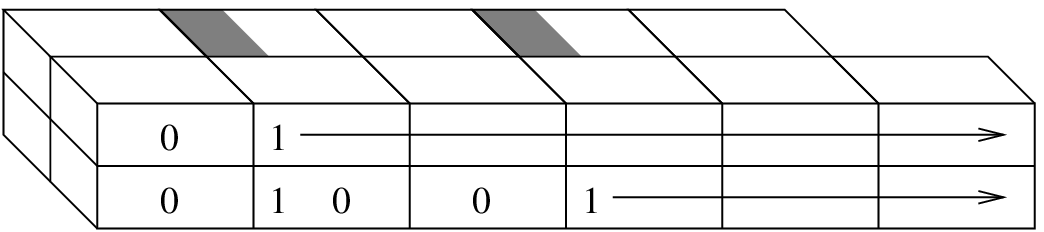}
\letteredcaption{b}{$\mu=(3,2)$, Level 1}\label{F-mp32-arblvl1}
\vspace{.5cm}
\psfig{figure=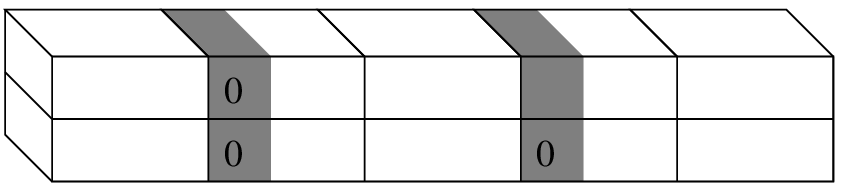}
\letteredcaption{b}{$\mu=(3,2)$, Level 2}\label{F-mp32-arblvl2}
\end{center}
\end{figure}
First we solve the 0th level as in Figure~\ref{F-mp32-arblvl0}, and
remove the two pivotal columns.  Now, in the second level, all  
 column operations are allowed, and row multiplication is allowed,
because the damage it does to the pivots can be repaired by a column
multiplication.  However, the rows cannot be added 
to each other.  This level is
solved as in Figure~\ref{F-mp32-arblvl1}.  For level 2 we get
the same row and column operations as in level 1, except that the
shaded columns  cannot  be 
added to other columns.  However this does not cause a problem,
since all but one entry in these columns has  already been killed by a
pivot on level 1.
\end{enumerate}
The final claim follows because all partitions of a number less than 6
are of one of these three forms.
\end{proof}

\begin{corollary}\label{C-ccft}
Computing the conjugacy classes in $P^{(m,n)}$ over a perfect field
reduces to matrix problems of finite type if, and only if, either $m<6$ or
$n<6$. In particular, computing conjugacy classes in the maximal
parabolics of the general linear group over a perfect field reduces to
finite type problems if, and only if, the dimension is less than 12.
\end{corollary}
\begin{proof}
The previous theorem, together with the results of
Sections~\ref{S-cent} to \ref{S-cocent}, show that we get finite type
problems if 
$m<6$.  By symmetry, this is also true for $n<6$.  By
Theorem~\ref{T-it}, $P^{(6,6)}$ involves a problem of infinite type and
it follows easily that $P^{(m,n)}$ does whenever $m,n\ge6$.
\end{proof}
When our  field is finite we get the following result.
\begin{corollary}\label{C-cnt}
Suppose that $k$ is finite of size $q$ and either $m<6$ or $n<6$.  Then the
number of conjugacy classes, and therefore the number of irreducible
characters, of $P^{(m,n)}$ is a polynomial in $q$ with  
integral coefficients.
\end{corollary}
\begin{proof} 
The number of solutions of the relevant matrix problems is
independent of $q$.  Hence this result follows immediately from the
well known fact that the number of characters of $\gl{n}q$ is a
polynomial in $q$ with  integral coefficients.
\end{proof}

Note that the proof of the  Theorem~\ref{T-ft} also provides a procedure
for solving these finite type problems, so this section gives an implicit
description of all conjugacy classes in the parabolic subgroups
mentioned in
Corollary~\ref{C-ccft}. 
For example, suppose the original eigenvalue is $\al$,  $\mu=\nu=(4,2)$,
and our orbit representative in $\mpr{k}{x}{\mu\times\nu}$ is 
$$\bmat \be x^2 & x \\ x & 1 \emat.$$ 
Then our conjugacy class representative in $P^{(6,6)}$ is
$$\left(\begin{array}{cccccc|cccccc}
\al&    &   &   &   &   &   &   &   &   &   &   \\
1  &\al &   &   &   &   &   &   &   &   &   &   \\
   &1   &\al&   &   &   &\be&   &   &   &   &   \\
   &    &1  &\al&   &   &   &\be&   &   &1  &   \\
   &    &   &   &\al&   &   &   &   &   &1  &   \\
   &    &   &   &1  &\al&1  &   &   &   &   &1  \\\hline

   &    &   &   &   &   &\al&   &   &   &   &   \\
   &    &   &   &   &   &1  &\al&   &   &   &   \\
   &    &   &   &   &   &   &1  &\al&   &   &   \\
   &    &   &   &   &   &   &   &1  &\al&   &   \\
   &    &   &   &   &   &   &   &   &   &\al&   \\
   &    &   &   &   &   &   &   &   &   &1  &\al
\end{array}\right),$$
where blank entries are zero.

\section{The affine general linear group}\label{S-ccagl}
We apply the results of the previous section to the affine general
linear groups \cite{mr47:6834,mr46:5475}.  The representation theory of these well known groups is
computed in \cite{MR83k:20017}. The affine general linear group of degree $n$
over the field $k$, $\aglnk$, is the semidirect product of $\glnk$ and
the row space $(k^n)^t$.  It can be realized as the subgroup  
$$\left(\begin{matrix} 1 & (k^n)^t\\ 0 &\glnk \end{matrix}\right)$$
of $\gl{n+1}{k}$.
The generalized Jordan normal form provides a set of conjugacy
class representatives 
for $\glnk$. Each can be written in  the form $N\ds E$
where $N$ has no eigenvalues equal to 1 and $E$ has eigenvalue 1.  There is a partition
$\la=(r^{l_r},\dots,2^{l_2},1^{l_1})$ so that
  \[E = J_\la(1) = \bigds_{i=1}^{r} E_i\]
where $E_i=J_i(1)^{\ds l_i}$.
\begin{theorem}\label{T-ccagl}
A set of conjugacy class representatives for $\aglnk$ is given by the matrices 
\begin{equation*}
\left(\begin{array}{c|cc}
        1&0&0\\ \hline
         &N&0\\
         &0&E\\
\end{array} \right) \mbox{ and }
\left(\begin{array}{c|cccccc}
        1&0&0     &\cdots &e       &\cdots&0\\ \hline
         &N&0     &\cdots &0       &\cdots&0\\
         &0&E_1   &       &        &      &0\\
         &\vdots& &\ddots &        &      &\vdots\\
         &0&      &       &E_i     &      &0\\
         &\vdots& &       &        &\ddots&\vdots\\
         &0&0     &\cdots &0       &\cdots&E_m \end{array} \right),
\end{equation*}
where $e = (1,0,0,\dots)$.
\end{theorem}

\begin{proof}
Let $Z=kI_{n+1}$ be the center of
$\gl{n+1}{k}$. Then $P^{(1,n)}=Z\cdot\aglnk$ and so the
conjugacy classes in $\aglnk$ are just the noncentral conjugacy classes in
$P^{(1,n)}$ intersected with $\aglnk$.  Hence we need the matrices
given by the proof of Theorem~\ref{T-ft} with $\mu=(1)$ and 1 in the
first summand of $L=\gl{1}{k}\ds\gl{n}{k}$.  The 
result is now immediate. 
\end{proof}

Let $k$ be a finite field.  We denote by $c_n$ the number of conjugacy
classes in $\glnk$ and use the convention that $\gl{0}{k}$ is the trivial
group.  For $d=0, 1,\dots,n$, we consider the conjugacy class
representatives of $\aglnk$ with an $e$ above a Jordan block of size
$d$.  Then $A = N\ds E$ is an
arbitrary conjugacy class representative of $\glnk$, except that it must have a
least one Jordan block of size $d$ and eigenvalue 1.  If you remove
one such block of  size $d$ from $A$, you get an arbitrary conjugacy class
representative of $\gl{n-d}{k}$, of which there are $c_{n-d}$. So the total
number of conjugacy class representatives of $\aglnk$ is 
\[\sum_{d=0}^n c_{n-d} = c_n + c_{n-1} +\dots + c_0.\] 
This agrees with the count of the number of irreducible 
characters of $\aglnk$ gotten by Zelevinsky \cite{MR83k:20017}.

\begin{acknowledgment}
I would like to thank my advisor, Jon Alperin.   I am also
grateful to Laszlo Babai, Cedric Bonnaf\'e, Joseph Chuang, Matt Frank, 
George Glauberman,  David Hemmer,  Paul Li, Paul Sally, and Larry
Wilson  for their suggestions and support.  This research was done
under the support of a Fulbright Graduate Fellowship.

\end{acknowledgment}


\begin{references}
\expandafter\ifx\csname natexlab\endcsname\relax\def\natexlab#1{#1}\fi

\bibitem{mr96m:20001}
J.~L. Alperin and Rowen~B. Bell.
\newblock {\em Groups and representations}.
\newblock Springer-Verlag, New York, 1995.

\bibitem{mr90k:20001}
Charles~W. Curtis and Irving Reiner.
\newblock {\em Methods of representation theory. {V}ol. {I}}.
\newblock John Wiley \& Sons Inc., New York, 1990.

\bibitem{mr52:14076}
P.~Deligne and G.~Lusztig.
\newblock Representations of reductive groups over finite fields.
\newblock {\em Ann. of Math. (2)}, 103(1):103--161, 1976.

\bibitem{mr94d:22016}
Yu.~A. Drozd.
\newblock Matrix problems, small reduction and representations of a class of
  mixed {L}ie groups.
\newblock In {\em Representations of algebras and related topics (Kyoto,
  1990)}, pages 225--249. Cambridge Univ. Press, Cambridge, 1992.

\bibitem{mr98e:16014}
P.~Gabriel and A.~V. Roiter.
\newblock {\em Representations of finite-dimensional algebras}.
\newblock Springer-Verlag, Berlin, 1997.
\newblock Translated from the Russian, With a chapter by B. Keller, Reprint of
  the 1992 English translation.

\bibitem{MR1669178}
L.~Hille and G.~R{\"o}hrle.
\newblock A classification of parabolic subgroups of classical groups with a
  finite number of orbits on the unipotent radical.
\newblock {\em Transform. Groups}, 4(1):35--52, 1999.

\bibitem{MR1692308}
Lutz Hille and Gerhard R{\"o}hrle.
\newblock On parabolic subgroups of classical groups with a finite number of
  orbits on the unipotent radial.
\newblock {\em C. R. Acad. Sci. Paris S\'er. I Math.}, 325(5):465--470, 1997.

\bibitem{MR99b:20011}
I.~M. Isaacs and Dikran Karagueuzian.
\newblock Conjugacy in groups of upper triangular matrices.
\newblock {\em J. Algebra}, 202(2):704--711, 1998.

\bibitem{MR1665003}
U.~J{\"u}rgens and G.~R{\"o}hrle.
\newblock Algorithmic modality analysis for parabolic groups.
\newblock {\em Geom. Dedicata}, 73(3):317--337, 1998.

\bibitem{MR86j:00003}
Serge Lang.
\newblock {\em Algebra}.
\newblock Addison-Wesley Publishing Co., Reading, Mass., third edition, 1993.

\bibitem{mr29:3477}
A.~I. Mal{'}cev.
\newblock {\em Foundations of linear algebra}.
\newblock W. H. Freeman \& Co., San Francisco, Calif.-London, 1963.

\bibitem{mr49:4877}
L.~A. Nazarova and A.~V. Ro{\u\i}ter.
\newblock Representations of partially ordered sets.
\newblock {\em Zap. Nau\v cn. Sem. Leningrad. Otdel. Mat. Inst. Steklov.
  (LOMI)}, 28:5--31, 1972.

\bibitem{MR54:360}
L.~A. {N}azarova and A.~V. {R}o{\u\i}ter.
\newblock {\em { {K}ategornye matrichnye zadachi i problema
  {B}rau\`era-{T}r\`ella.}}
\newblock Izdat. ``Naukova Dumka'', Kiev, 1973.

\bibitem{MR99f:14063}
Vladimir Popov and Gerhard R{\"o}hrle.
\newblock On the number of orbits of a parabolic subgroup on its unipotent
  radical.
\newblock In {\em Algebraic groups and Lie groups}, pages 297--320. Cambridge
  Univ. Press, Cambridge, 1997.

\bibitem{MR1617826}
Vladimir~L. Popov.
\newblock A finiteness theorem for parabolic subgroups of fixed modality.
\newblock {\em Indag. Math. (N.S.)}, 8(1):125--132, 1997.

\bibitem{mr93j:20092}
Roger Richardson, Gerhard R{\"o}hrle, and Robert Steinberg.
\newblock Parabolic subgroups with abelian unipotent radical.
\newblock {\em Invent. Math.}, 110(3):649--671, 1992.

\bibitem{MR97c:20070}
Gerhard R{\"o}hrle.
\newblock Parabolic subgroups of positive modality.
\newblock {\em Geom. Dedicata}, 60(2):163--186, 1996.

\bibitem{MR99b:20081}
Gerhard R{\"o}hrle.
\newblock Maximal parabolic subgroups in classical groups are of modality zero.
\newblock {\em Geom. Dedicata}, 66(1):51--64, 1997.

\bibitem{MR99e:20060}
Gerhard R{\"o}hrle.
\newblock A note on the modality of parabolic subgroups.
\newblock {\em Indag. Math. (N.S.)}, 8(4):549--559, 1997.

\bibitem{MR1669615}
Gerhard R{\"o}hrle.
\newblock On the modality of parabolic subgroups of linear algebraic groups.
\newblock {\em Manuscripta Math.}, 98(1):9--20, 1999.

\bibitem{mr47:6834}
Louis Solomon.
\newblock On the affine group over a finite field.
\newblock In {\em Representation theory of finite groups and related topics
  (Proc. Sympos. Pure Math., Vol. XXI, Univ. Wisconsin, Madison, Wis., 1970)},
  pages 145--147. Amer. Math. Soc., Providence, R.I., 1971.

\bibitem{mr46:5475}
Louis Solomon.
\newblock The affine group. {I}. {B}ruhat decomposition.
\newblock {\em J. Algebra}, 20:512--539, 1972.

\bibitem{MR42:3091}
T.~A. Springer and R.~Steinberg.
\newblock Conjugacy classes.
\newblock In {\em Seminar on Algebraic Groups and Related Finite Groups (The
  Institute for Advanced Study, Princeton, N.J., 1968/69)}, pages 167--266.
  Springer, Berlin, 1970.
\newblock Lecture Notes in Mathematics, Vol. 131.

\bibitem{MR83k:20017}
Andrey~V. Zelevinsky.
\newblock {\em Representations of finite classical groups}.
\newblock Springer-Verlag, Berlin, 1981.
\newblock A Hopf algebra approach.

\end{references}

\end{article}

\end{document}